\documentclass[11pt]{article}
\usepackage{graphicx}
\usepackage{amsmath}
\usepackage{amssymb}
\usepackage{amsthm}
\usepackage{appendix}
\usepackage{url}
\newcommand\RR{\mathbb{R}}

\newcommand\al\alpha
\newcommand\be\beta
\newcommand\de\delta
\newcommand\ep\varepsilon
\newcommand\tha\theta
\newcommand\ka\kappa
\newcommand\la\lambda
\newcommand\om\omega

\newcommand\iy\infty
\newcommand\pa\partial

\newcommand{\hyp}[5]{\,\mbox{}_{#1}F_{#2}\!\left(\genfrac{}{}{0pt}{}{#3}{#4};#5\right)}

\numberwithin{equation}{section}
\newtheorem{theorem}{Theorem}[section]

\newtheorem{Remark}[theorem]{Remark}

\begin{document}

\title{The orthogonal Shmaliy polynomials are Hahn polynomials.}
\author{Enno Diekema \footnote{email adress: e.diekema@gmail.com}}
\maketitle

\begin{abstract}
\noindent
Morales-Mendoza et al. present in 2013 a new class of discrete orthogonal polynomials. They use these polynomials to design an unbiased FIR filter. In their paper they make the statement that a representation of the polynomials via hypergeometric functions is unknown. However Shakibaei Asli et al. found in 2017 a $_3F_2$ hypergeometric function representation. In this paper it is shown that the "new" orthogonal  polynomials belong to the class of the Hahn polynomials. The transfer function of the unbiased FIR filter has been determined. In the second part, an orthogonal unbiased FIR filter is designed following the method of the orthogonal derivative described in the thesis of the author. In the Appendix a table with a number of possible transformations of the Hahn polynomials is given.
\end{abstract}

\

\section{Introduction}
\setlength{\parindent}{0cm}

The author during his work at an hospital had to calculate the derivative of a given function concerning the pressure in the eye. The signal was seriously disturbed by noise so an averaging routine had to be used. He used the well-known least square procedure. The result was a formula for a discrete low pass filter based on a summation with discrete orthogonal polynomials. In this way the $\text{\em orthogonal difference}$ arose. He then expands the theory to the continuous case and so he found a formula for the $\text{\em orthogonal derivative}$ based on an integral with continuous orthogonal polynomials \cite{3}, \cite{4}. 

The formula found for the $n-$th orthogonal derivative and $n-$th orthogonal difference are both dependent of the order of the used orthogonal polynomials. When the order is set to zero the filter becomes a low pass filter based on orthogonal polynomials. In this paper we consider only the discrete filter. 

In \cite[section 4]{4} an discrete orthogonal filter is treated and in \cite[section 5]{4} the theory is applied using the Hahn polynomials. The filter is a finite impulse response (FIR) filter. The author was unaware of the papers of Shmaliy in which Shmaliy et al. developed a theory for an unbiased FIR filter (UFIR) with an orthogonal polynomial \cite{5}. In \cite{6} the authors claim that the orthogonal polynomials found are new orthogonal polynomials. In \cite{7} the authors found that the "new" orthogonal polynomials can be described as a hypergeometric function.

In section 2 it is shown that the UFIR functions derived by Shmaliy are Hahn polynomials. In section 3 the transfer function of the UFIR is derived. In section 4 the orthogonal UFIR filter based on the theory of the orthogonal derivative is treated.

 \newpage

\section{The UFIR function as a Hahn polynomial}
The $n-$th degree discrete Shmaliy polynomials are defined as
\[
h_m(n,N)=\sum_{i=0}^m a_{im}(N)n^i
\]
with $0 \leq n \leq N-1$. The coefficients $a_{im}$ are defined as
\[
a_{im}(N)=(-1)^i\dfrac{M^{(m)}_{(i+1),1}(N)}{|H_m(N)|}
\]
where $H_m(N)|$ is the determinant and $M^{(m)}_{(i+1),1}(N)$ is the minor of the Hankel matrix $H_m(N)$
\[
H_m(N)=\left(\begin{array}{cccc}
	c_0&c_1&\dots&c_m \\
	c_1&c_2&\dots&c_{m+1} \\
	\dots&\dots&\dots&\dots \\
	c_m&c_{m+1}&\dots&c_{2m}
\end{array}
\right)
\]
The elements of the Hankel matrix are given by the power series
\[
c_k(N)=\sum_{i=0}^{N-1}i^k=\dfrac{1}{k+1}\left(B_{k+1}(N)-B_{k+1}\right)
\]
where $B_k(x)$ is the Bernoulli polynomial and $B_k=B_k(0)$ is the Bernoulli number. For the Shmaliy polynomial we get
\[
h_m(n,N)=\sum_{i=0}^m(-1)^i\dfrac{M^{(m)}_{(i+1),1}(N)}{|H_m(N)|}
\]
In \cite{6} the following properties of the discrete polynomials are proved
\begin{align*}
&\sum_{n=0}^{N-1}h_m(n,N)=1 \\
&\sum_{n=0}^{N-1}x^q h_m(n,N)=0 \qquad 1 \leq q \leq n \\
&\sum_{n=0}^{N-1}h_m^2(n,N)<\infty
\end{align*}
The orthogonality property is
\[
\sum_{n=0}^{N-1}\rho(n,N)h_m(n,N)h_q(n,N)=(d_m)^2\delta_{mq}
\]
with weight function
\[
\rho(n,N)=\dfrac{2n}{N(N-1)}
\]
and 
\[
(d_m)^2=\dfrac{(m+1)(N-m-1)}{N(N)_{m+1}}=\dfrac{(m+1)\Gamma(N-1)\Gamma(N)}{\Gamma(N-m-1)N\Gamma(N+m+1)}
\]
They found the recurrence relation
\[
h_m(n,N)=2\dfrac{m^2(2N-1)-(4m^2-1)n}{m(2m-1)(N+m)}h_{m-1}(n,N)-
\dfrac{(2m+1)(N-m)}{(2m-1)(N+m)}h_{m-1}(n,N)
\]
with $m \geq1,\ h_{-1}(n,N)=0$ and $h_0(n,N)=\dfrac{1}{N}$. The first four functions are
\begin{align*}
&h_0(n,N)=\dfrac{1}{N} \\
&h_1(n,N)=\dfrac{2(2N-1)-6n}{N(N+1)} \\
&h_2(n,N)=\dfrac{3(3N^2-3N+2)-18(2N-1)n+30n^2}{N(N+1)(N+2)} \\
&h_3(n,N)=\dfrac{120(2N-1)n^2-140n^3}{N(N+1)(N+2)(N+3)}-
\dfrac{8(2N^3-3N^2+7N-3)+20(6N^2-6N+5)n}{N(N+1)(N+2)(N+3)}
\end{align*}
In \cite{7} the authors prove that $h_m(n,N)$ can be written as a hypergeometric function
\begin{equation}
h_m(n,N)=(-1)^m\dfrac{(m+1)(n-m)_m(N-n)_m}{m!(N)_{m+1}}
\hyp32{-m,n+1,1-N+n}{n-m,1-N-m+n}{1}
\label{2.1}
\end{equation}
We found by trial and error the following much simpler hypergeometric function
\begin{equation}
h_m(n,N)=\dfrac{(m+1)^2}{N}\hyp32{-m,n+1,m+2}{2,1+N}{1}
\label{2.2}
\end{equation}
To prove that \eqref{2.1} is equal \eqref{2.2} we use the following properties of the Pochhammer symbols and the Gamma functions
\[
(-1-m)_m=\frac{\Gamma(-1)}{\Gamma(-1-m)}=(-1)^m\Gamma(m+2)
\]
\[
\dfrac{\Gamma(N-n+m)\Gamma(1-N-m+n)}{\Gamma(N-n)\Gamma(1-N+n)}=
\dfrac{\sin(N-n)\pi}{\sin(N-n+m)\pi}=\dfrac{1}{\cos m\pi}=(-1)^m
\]
\[
\dfrac{\Gamma(N+1)\Gamma(-N)}{\Gamma(1+N+m)\Gamma(-N-m)}=(-1)^m
\]
Application of \eqref{A1} to \eqref{2.1} gives
\begin{align*}
h_m(n,N)&=(-1)^m\dfrac{(m+1)(n-m)_m(N-n)_m}{m!(N)_{m+1}}
\dfrac{(-1-m)_m(N-m)_m}{(n-m)_m(1-N-m+n)_m} \\
&\qquad\qquad\qquad\qquad\qquad\qquad\qquad\qquad\qquad\qquad\qquad\qquad \hyp32{-m,n+1,m+2}{2,1+N}{1} \\
&=\dfrac{(-1)^m(m+1)}{\Gamma(m+1)}(-1-m)_m\dfrac{(N-n)_m}{(1-N-m+n)_m}
\dfrac{(-N-m)_m}{(N)_{m+1}} \\
&\qquad\qquad\qquad\qquad\qquad\qquad\qquad\qquad\qquad\qquad\qquad\qquad
\hyp32{-m,n+1,m+2}{2,1+N}{1}
\end{align*}
\begin{align*}
\qquad\ \ \ &=\left[\dfrac{(-1)^m(m+1)}{\Gamma(m+1)}(-1-m)_m\right]
\left[\dfrac{\Gamma(N-n+m)\Gamma(1-N-m+n)}{\Gamma(N-n)\Gamma(1-N+n)}\right]
\left[\dfrac{\Gamma(N)}{\Gamma(N+1)}\right] \\
&\qquad\qquad\qquad\qquad\qquad\qquad\ \ \left[\dfrac{\Gamma(N+1)\Gamma(-N)}{\Gamma(1+N+m)\Gamma(-N-m)}\right]
\hyp32{-m,n+1,m+2}{2,1+N}{1} \\
&=\Big[(m+1)^2\Big]\Big[(-1)^m\Big]\Big[\dfrac{1}{N}\Big]\Big[(-1)^m\Big]\hyp32{-m,n+1,m+2}{2,1+N}{1} \\
&=\dfrac{(m+1)^2}{N}\hyp32{-m,n+1,m+2}{2,1+N}{1}
\end{align*}
and this is \eqref{2.2}. This is one of the Thomae transformations for the $_3F_2$  hypergeometric function with unit argument. In Appendix A the standard transformations are given. 

Both hypergeometric functions can be written as Hahn polynomials. \cite[(9.5.1)]{9} gives as definition
\[
Q_m(y;\alpha,\beta,M)=\hyp32{-m.m+\alpha+\beta+1,-y}{\alpha+1,-M}{1} \qquad m=0,1,2,\dots,M
\]
which are orthogonal as $\alpha>-1$ and $\beta>-1$, or $\alpha<-M$ and $\beta<-M$. These are polynomials of $m$-th degree.

Setting $\alpha=1$, $\beta=0$, $M=-1-N$ and $y=-n-1$ gives
\[
h_m(n,N)=\dfrac{(m+1)^2}{N}Q_m(-1-n;1,0,-1-N)
\]
With the transformations of the $_3F_2$ hypergeometric function there are 31 transformations of the Hahn polynomials. See Appendix B. We choose
\[
Q_m(y;\alpha,\beta,N)=\dfrac{(-m-\beta)_m(-N-m-\alpha-\beta -1)_m}{(\alpha +1)_m(-N)_m} Q_m(y-N-\beta ;\beta ,\alpha ,-N-\alpha-\beta -2)
\]
Application gives
\[
h_m(n,N)=\dfrac{(m+1)^2}{N}\dfrac{(-m)_m(N-1-m)_m}{(2)_m(N+1)_m}Q_m(N-n;0,1,N-2)
\]
Simplifying with $(-m)_m=(-1)^m\Gamma(m+1)$ gives at last
\[
h_m(n,N)=\dfrac{(m+1)(N+m+1)}{N(N+1)}Q_m(N-n;0,1,N-2)
\]
The conclusion is that the function $h_m(n,N)$ is a well known orthogonal Hahn polynomial.

\section{The general transfer function of the UFIR Shmaliy function}
In \cite{10} Shmialy et al. compute the transfer functions of the Unbiased FIR function $h_m(N,N)$ for the values $m=1,\ 2$ and $3$. In this section a general formula for the transfer function is derived. We start with \eqref{2.2}
\[
h_m(n,N)=\dfrac{(m+1)^2}{N}\hyp32{-m,n+1,m+2}{2,N+1}{1}
\]
and write the hypergeometric function as a summation
\[
h_m(n,N)=\dfrac{(m+1)^2}{N}\sum_{k=0}^m\dfrac{(-m)_k(n+1)_k(m+2)_k}{(2)_k(N+1)_k}\dfrac{1}{k!}
\]
For computing the transfer function we use the $z$-transform with $z=\exp(j\omega T)$
\begin{align*}
H_m(z,N)&=\sum_{n=0}^{N-1}h_m(n,N)z^{-n} \\
&=\sum_{n=0}^{N-1}\dfrac{(m+1)^2}{N}\sum_{k=0}^m\dfrac{(-m)_k(n+1)_k(m+2)_k}{(2)_k(N+1)_k}\dfrac{1}{k!}z^{-n}
\end{align*}
After using $(1+n)_k=(1+k)_n(1)_k/(1)_n$ and interchanging the summations we get
\[
H_m(z,N)=\dfrac{(m+1)^2}{N}\sum_{k=0}^m\dfrac{(-m)_k(m+2)_k(1)_k}{(2)_k(N+1)_k}\dfrac{1}{k!}\sum_{n=0}^{N-1}\dfrac{(1+k)_n}{(1)_n}z^{-n}
\]
The last summation is known. 
\begin{multline*}
H_m(z,N)=\dfrac{(m+1)^2}{N}\left(\dfrac{z}{z-1}\right)\sum_{k=0}^m\dfrac{(-m)_k(m+2)_k(1)_k}{(2)_k(N+1)_k}\dfrac{1}{k!}\left(\dfrac{z}{z-1}\right)^k \\
\left[1-z^{-N}(N+1)_k\dfrac{1}{k!}\hyp21{-k,N}{N+1}{\dfrac{1}{z}}\right]
\end{multline*}
Applying a Gauss transformation to the hypergeometric function gives
\begin{multline*}
H_m(z,N)=\dfrac{(m+1)^2}{N}\left(\dfrac{z}{z-1}\right)\sum_{k=0}^m\dfrac{(-m)_k(m+2)_k(1)_k}{(2)_k(N+1)_k}\dfrac{1}{k!}\left(\dfrac{z}{z-1}\right)^k- \\
-z^{-N}\dfrac{(m+1)^2}{N}\sum_{k=0}^m\dfrac{(-m)_k(m+2)_k}{(2)_k}\dfrac{1}{k!}
\sum_{i=0}^\infty\dfrac{(N+1+k)_i(1)_i}{(N+1)_i}\dfrac{1}{i!}\left(\dfrac{1}{z}\right)^i
\end{multline*}
After using $(N+1+k)_i=(N+1+i)_k(N+1)_i/(N+1)_k$, \ interchanging the summations and writing the summation over $k$ as a hypergeometric function we get 
\begin{multline}
H_m(z,N)=\dfrac{(m+1)^2}{N}\left(\dfrac{z}{z-1}\right)\sum_{k=0}^m\dfrac{(-m)_k(m+2)_k(1)_k}{(2)_k(N+1)_k}\dfrac{1}{k!}\left(\dfrac{z}{z-1}\right)^k- \\
-z^{-N}\dfrac{(m+1)^2}{N}\sum_{i=0}^\infty\left(\dfrac{1}{z}\right)^i
\hyp32{-m,m+2,N+1+i}{2,N+1}{1}
\label{3.1}
\end{multline}
For the hypergeometric function we use \eqref{A6}.
\begin{multline*}
\hyp32{-m,m+2,N+1+i}{2,N+1}{1}=\dfrac{(-m)_m}{(2)_m}\hyp32{-i,m,m+2}{1,N+1}{1}= \\
=\dfrac{(-1)^m}{(m+1)}\hyp32{-i,m,m+2}{1,N+1}{1}=
\dfrac{(-1)^i}{(m+1)}\sum_{k=0}^\infty\dfrac{(-i)_k(-m)_k(m+2)_k}{(1)_k(N+1)_k}\dfrac{1}{k!}
\end{multline*}
Using $(-i)_k=(-1)^k\dfrac{(1)_i}{\Gamma(1-k)(1-k)_i}$ gives
\[
\hyp32{-m,m+2,N+1+i}{2,N+1}{1}=\dfrac{(-1)^m(-1)_i}{(m+1)}
\sum_{k=0}^\infty\dfrac{(-m)_k(m+2)_k}{(1)_k(N+1)_k\Gamma(1-k)(1-k)_i}\dfrac{1}{k!}
(-1)^k
\]
Substitution in \eqref{3.1} and interchanging the summation in the last hypergeometric function gives
\begin{multline}
H_m(z,N)=\dfrac{(m+1)^2}{N}\left(\dfrac{z}{z-1}\right)\sum_{k=0}^m\dfrac{(-m)_k(m+2)_k(1)_k}{(2)_k(N+1)_k}\dfrac{1}{k!}\left(\dfrac{z}{z-1}\right)^k- \\
-z^{-N}\dfrac{(m+1)^2}{N}\sum_{k=0}^m\dfrac{(-m)_k(m+2)_k}{(1)_k(N+1)_k}\dfrac{1}{k!}
\dfrac{(-1)^k}{\Gamma(1-k)}
\sum_{i=0}^\infty \dfrac{(1)_i(1)_i}{(1-k)_i}\left(\dfrac{1}{z}\right)^i
\label{3.2}
\end{multline}
For the last summation we can use the following property of the hypergeometric function 
\begin{multline*}
\dfrac{1}{\Gamma(-M)} \ _{p+1}F_p
\left(\begin{array}{l}
	a_0,\dots,a_p \\
	-M,b_2,\dots,b_p
\end{array};x\right)= \\
=\dfrac{x^{M+1}(a_0)_{M+1}\dots(a_p)_{M+1}}{\Gamma(M+2)(b_2)_{M+1}\dots(b_p)_{M+1}}
\ _{p+1}F_p
\left(\begin{array}{l}
	a_0+M+1,\dots,a_p+M+1 \\
	M+2,b_2+M+1,\dots,b_p+M+1
\end{array};x\right)
\end{multline*}
with $M$ a positive integer. For the proof of this formula see \cite[Lemma 2]{11}. The result is
\begin{multline*}
\dfrac{(-1)^k}{\Gamma(1-k)}
\sum_{i=0}^\infty \dfrac{(1)_i(1)_i}{(1-k)_i}\left(\dfrac{1}{z}\right)^i=
(-1)^k\left(\dfrac{1}{z}\right)^k(1)_k\sum_{i=0}^\infty(1+k)_i\dfrac{1}{i!}
\left(\dfrac{1}{z}\right)^i= \\
=(-1)^k\left(\dfrac{1}{z}\right)^k(1)_k\left(\dfrac{z}{z-1}\right)^{k+1}=
\left(\dfrac{z}{z-1}\right)(1)_k\left(\dfrac{1}{1-z}\right)^k
\end{multline*}
Substitution in \eqref{3.2} gives
\begin{multline*}
H_m(z,N)=\dfrac{(m+1)^2}{N}\left(\dfrac{z}{z-1}\right)
\sum_{k=0}^m\dfrac{(-m)_k(m+2)_k(1)_k}{(2)_k(N+1)_k}\dfrac{1}{k!}\left(\dfrac{z}{z-1}\right)^k- \\
-z^{-N}\dfrac{(m+1)}{N}(-1)^m\left(\dfrac{z}{z-1}\right)
\sum_{k=0}^m\dfrac{(-m)_k(m+2)_k}{(N+1)_k}\dfrac{1}{k!}\left(\dfrac{1}{1-z}\right)^k
\end{multline*}
The first summation is known.
\[
\sum_{k=0}^m\dfrac{(-m)_k(m+2)_k(1)_k}{(2)_k(N+1)_k}\dfrac{1}{k!}\left(\dfrac{z}{z-1}\right)^k=\dfrac{N}{(1+m)^2}\left(\dfrac{z-1}{z}\right)\left[1-\hyp21{-1-m,1+m}{N}{\dfrac{z}{z-1}}\right]
\]
Substitution gives for the transfer function
\[
H_m(z,N)=1-\hyp21{-1-m,1+m}{N}{\dfrac{z}{z-1}}-(-1)^m\dfrac{(m+1)}{N}\dfrac{z^{1-N}}{z-1}
\hyp21{-m,m+2}{N+1}{\dfrac{1}{1-z}}
\]
Application of Gauss transformations to both hypergeometric functions results in
\begin{multline*}
H_m(z,N)=1-\left(\dfrac{1}{1-z}\right)^{m+1}\hyp21{-1-m,N-m-1}{N}{z}+ \\
+\dfrac{(m+1)}{N}z^{-N}\left(\dfrac{z }{1-z}\right)^{m+1}\hyp21{-m,N-m-1}{N+1}{\dfrac{1}{z}}
\end{multline*}
Because $|z| \leq 1$ we use for the second hypergeometric function the transformation
\[
\hyp21{-m,N-m-1}{N+1}{\dfrac{1}{z}}
=\dfrac{\Gamma(N+1)\Gamma(N-1)}{\Gamma(N-m-1)\Gamma(N+m+1)}\left(-\dfrac{1}{z}\right)^m
\hyp21{-m,-N-m}{2-N}{z}
\]
Substitution gives
\begin{multline*}
H_m(z,N)=1-\left(\dfrac{1}{1-z}\right)^{m+1}\hyp21{-1-m,N-m-1}{N}{z}+ \\
+z^{1-N}\left(\dfrac{1}{1-z}\right)^{m+1}
\dfrac{(m+1)\Gamma(N)\Gamma(N-1)}{\Gamma(N-m-1)\Gamma(N+m+1)}(-1)^m
\hyp21{-m,-N-m}{2-N}{z}
 \end{multline*}
Both hypergeometric functions can be written as Jacobi polynomials. There are a lot of possibilities. We take the following transformations
\begin{align*}
&\hyp21{-1-m,N-m-1}{N}{z}
=\dfrac{\Gamma(N)\Gamma(m+2)}{\Gamma(N+m+1)}P_{m+1}^{(N-1,-2-2m)}(1-2z) \\
&\hyp21{-m,-N-m}{2-N}{z}
=(-1)^m\dfrac{\Gamma(N-m-1)\Gamma(m+1)}{\Gamma(N-1)}P_m^{(1-N,-2-2m)}(1-2z)
\end{align*}
Substitution gives at last
\begin{multline*}
H_m(z,N)=1-\left(\dfrac{1}{1-z}\right)^{m+1}\dfrac{\Gamma(N)\Gamma(m+2)}{\Gamma(N+m+1)} \\
\left[P_{m+1}^{(N-1,-2-2m)}(1-2z)-z^{1-N}P_m^{(1-N,-2-2m)}(1-2z)\right]
\end{multline*}
If we replace $z$ with $e^{j\omega T}$ we get the transfer function of the UFIR function.
\begin{multline*}
H_m(\omega\,T,N)=1-\left(\dfrac{1}{1-e^{j\omega T}}\right)^{m+1}
\dfrac{\Gamma(N)\Gamma(m+2)}{\Gamma(N+m+1)} \\
\left[P_{m+1}^{(N-1,-2-2m)}(1-2e^{j\omega T})-e^{j(1-N)\omega T}
P_m^{(1-N,-2-2m)}(1-2e^{j\omega T})\right]
\end{multline*}
The following figure shows the modulus of the transfer function for the values of $m=1,\, 2,\, 3$. See also \cite[Fig. 4]{12}.
\begin{figure}[ht]
\centering     
\includegraphics[width=8cm]{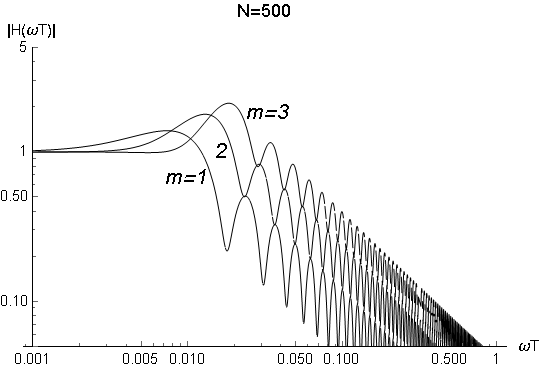} 
\caption{Absolute value of the transfer function for $N=500$ and for different values for $m$ as a function of $\omega\,T$.}
\label{Figure 1}
\end{figure}

It turns out that for small $\omega\, T$ a problem arises because both Jacobi polynomials  between the square brackets approach about the same value. Then there is a great loss of decimals for the difference of these terms. The result is unpredictable. See the left picture in Figure 2.

\begin{figure}[ht]
\centering     
\includegraphics[width=16cm]{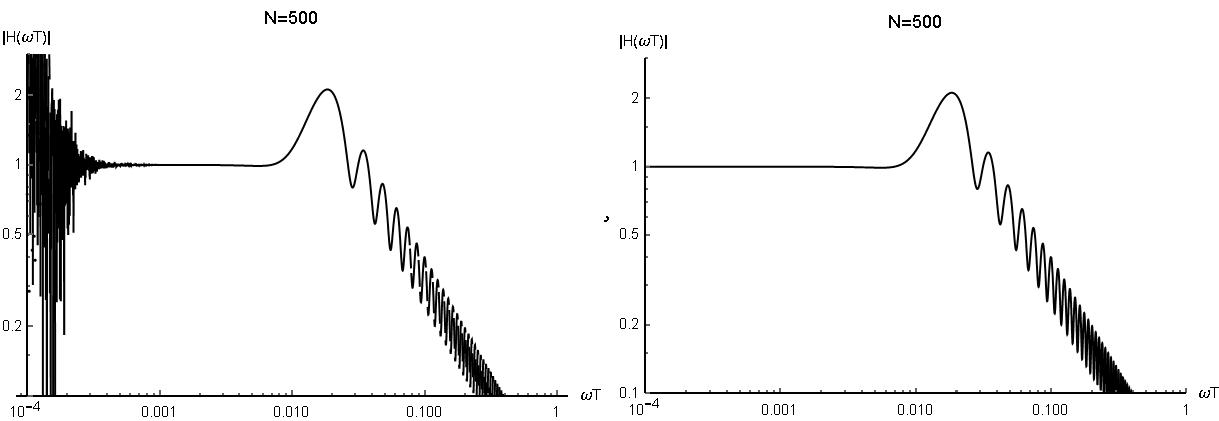} 
\caption{Absolute value of the transfer functions for $N=500$ and $m=3$. The horizontal axis begins at $\omega\, T=10^{-4}$. Left picture: 9 decimals are used. Right picture: 50 decimals are used.}
\label{Figure 2}
\end{figure}
Therefore, the calculation accuracy should be increased. The result is shown in Figure 2. In the left picture the standard accuracy of 9 decimal places is used. In the right picture the accuracy has been increased to 20 decimal places. This means that for signals with a low frequency content ($\omega\, T$ small), a greater number of decimals should be used for the filter to function properly. 

\section{The orthogonal unbiased FIR filter}
In \cite{3} and in his thesis \cite{4}, the author describes the theory of the orthogonal derivative. The theory is developed for both the continuous case and the discrete case. It is not discussed what happens when the order is set equal to zero. It will appear that a low-pass filter is then created. Since a UFIR involving Hahn polynomials has been discussed in the previous sections, we also use the Hahn polynomials to investigate the properties of the low-pass filter. See also \cite{14}.

\

For the definition of the orthogonal derivative we use \cite[Theorem 2.3.1]{4}.
\begin{theorem}
For some $n$ let $p_{n}$ be an orthogonal polynomial of degree $n$ with respect to the orthogonality measure $\mu $. Let $x\in \RR$. Let $I$  be a closed interval such that, for some $\varepsilon >0,$ $\
x+\delta v\in I$ if $0\leq \delta \leq \varepsilon $ and $v\in $ supp$(\mu) $. Let $f$ be a continuous function on $I$ such that its derivatives of order $1,2,\dots,n$ at $x$ exist. In addition, if $I$ is unbounded, assume that $f$ is of at most polynomial growth on $I$. Then
\[
f^{(n) }(x) =\lim\limits_{\delta\downarrow 0}D_{\delta }^{n}[f](x)
\]
where
\[
D_{\delta }^{n}[f](x) =\dfrac{k_{n}n!}{h_{n}}\dfrac{1}{\delta^{n}}\int_{\RR}
f(x+\delta \xi) p_{n}(\xi) d\mu(\xi)
\]
where the integral converges absolutely.
\end{theorem} 
In the rest of the paper the limit will be omitted. We set $\delta=1$. In that case the orthogonal derivative becomes an approximate orthogonal derivative. The theory for the discrete case using the Hahn polynomials is used with $\mu$ has a finite support ${x_0,x_1,x_2,\dots,x_N}$ with corresponding weights $w_0,w_1,w_2,\dots,w_N$. Then the definition becomes
\[
f^{(n)}(x)=\dfrac{k_nn!}{h_n}\sum_{k=0}^{N-1}f(x+k)Q_n(x;\alpha,\beta,N)w(x;\alpha,\beta,N)
\]
For the Hahn polynomials we use the definition in \cite[(9.5.1)]{9}
\[
Q_{n}(x;\alpha ,\beta ,N):=\hyp32{-x,-n,n+\alpha +\beta +1}{\alpha +1,-N}{1} \qquad n=0,1,2,\dots,N
\]
Here the hypergeometric series is assumed to be a sum from $k=0$ to $n$, by which the bottom parameter $-N$ will not cause a singularity. For $\alpha,\beta>-1$ or $\alpha,\beta<-N$ they are orthogonal on $0,1,2.\dots,N$ with respect to the weights
\[
w(x;\alpha ,\beta ,N) =\dfrac{(\alpha+1)_{x}(\beta +1)_{N-x}}{x!(N-x)!}
\]
and
\[
\dfrac{k_{n}}{h_{n}}=(-1) ^{n}\dfrac{(2n+\alpha +\beta +1) }{(\beta+1) _{n}}
\dfrac{(n+\alpha +\beta +1) _{n}}{(n+\alpha +\beta +1) _{N+1}}\dfrac{N!}{n!}
\]
For all formula in this section we have  $x=0,1,\dots,N$. So $x$ is an integer. Then we have $1/x!=0$ for $x=-1,-2,\dots,-n$ and $1/(N-n-x)!=0$ for $x=N,N-1,\dots,N-n+1$.

Setting $n=0$ gives
\begin{equation}
y(x,N)=\dfrac{\Gamma(N+1)}{(\alpha+2)_N}\sum_{m=0}^{N-1} f(x+m)
\dfrac{(\alpha+1)_m(\beta+1)_{N-m}}{(N-m)!m!}
\label{4.0}
\end{equation}
We need this equation for the $z$ transform. For simplicity we set $\beta=0$. The result is
\[
y(x,N)=\dfrac{\Gamma(N+1)}{(\alpha+2)_N}\sum_{m=0}^{N-1} f(x+m)\dfrac{(\alpha+1)_m}{m!}
\]
If $\alpha$ is an integer, the impulse response is a polynomial of order $\alpha$.
This means that the behaviour of this function can be compared with the function of Shmaliy.

\

For $\alpha=0$ the filter gives a simple average.
\[
y(x,N)=\dfrac{1}{N+1}\sum_{m=0}^{N-1} f(x+m)
\]
For $\alpha=1$ the filter becomes
\[
y(x,N)=\dfrac{2}{(N+1)(N+2)}\sum_{m=0}^{N-1} (m+1)f(x+m)
\]

\

For computing the transfer function we use the $z$-transform of \eqref{4.0}.
\begin{equation}
H(z,N)=\dfrac{\Gamma(N+1)}{(\alpha+2)_N}\sum_{m=0}^{N-1}
\dfrac{(\alpha+1)_m(\beta+1)_{N-m}}{(N-m)!m!}z^{-m}
\label{4.0a}
\end{equation}
The summation is well known.
\begin{equation}
H(z,N)=\dfrac{(\beta+1)_N}{(\alpha+2)_N}
\hyp21{-N,\alpha+1}{-\beta-N}{\dfrac{1}{z}}-\dfrac{(\alpha+1)_N}{(\alpha+2)_N}z^{-N}
\label{4.0b}
\end{equation}
For simplicity we want $\beta=0$. But for $\beta=0$ this formula is not valid. Setting $\beta=0$ in \eqref{4.0a} gives
\begin{equation}
H(z,N,\alpha)=\left(\dfrac{z}{z-1}\right)^{\alpha+1}\dfrac{\Gamma(N+1)}{(\alpha+2)_N}-
\dfrac{\alpha+1}{N+\alpha+1}z^{-N}\left(\dfrac{z}{z-1}\right)^{\alpha+1}
\hyp21{-\alpha,N}{N+1}{\dfrac{1}{z}}
\label{4.0c}
\end{equation}

\

Using \cite[2.10.(4)]{13} we get for the hypergeometric function
\begin{multline*}
\hyp21{-\alpha,N}{N+1}{\dfrac{1}{z}}=
\dfrac{\Gamma(N+1)\Gamma(\alpha+1)}{\Gamma(N+\alpha+1)}z^N- \\
-\dfrac{N}{(\alpha+1)}z^N(z-1)^{\alpha+1}
\hyp21{N+\alpha+1,\alpha+1}{\alpha+2}{1-z}
\end{multline*}
Substitution in \eqref{4.0c}gives
\[
H(z,N,\alpha)=\dfrac{N}{N+\alpha+1}z^{\alpha+1}\hyp21{N+\alpha+1,\alpha+1}{\alpha+2}{1-z}
\]
Application of a Gauss transformation results in a bounded hypergeoetric function
\begin{equation}
H(z,N,\alpha)=\dfrac{N}{N+\alpha+1}z^{1-N}\hyp21{1-N,1}{\alpha+2}{1-z}
\label{4.1}
\end{equation}
The hypergeometric function can also be written as a Jacobi polynomial
\[
H(z,N,\alpha)=\dfrac{\Gamma(N+1)}{(\alpha+2)_N}z^{1-N}P_{N-1}^{(\alpha+1,-\alpha-N)}(2z-1)
\]
If we replace $z$ with $e^{j\omega\, T}$ we get the transfer function of the filter.
\[
H(z,N,\alpha)=\dfrac{N}{N+\alpha+1}e^{j(1-N)\omega\, T}\hyp21{1-N,1}{\alpha+2}{1-e^{j\omega\, T}}
\]
Taking the limit for $\omega \rightarrow 0$ and for $\omega \rightarrow 2\pi$ we get
\begin{equation}
\lim_{\omega \rightarrow 0,2\pi}H(\omega,N)=\dfrac{N}{N+\alpha+1}
\label{4.2a}
\end{equation}
For large $N$ this limit goes to $1$. This means that the filter is a LP filter with parameter $\alpha$. In the Figure 3 we show this transfer function with $N=200$ and different values for $\alpha$. For high frequencies the filter behaves like a first order LP filter. 
\begin{figure}[ht]
\centering     
\includegraphics[width=10cm]{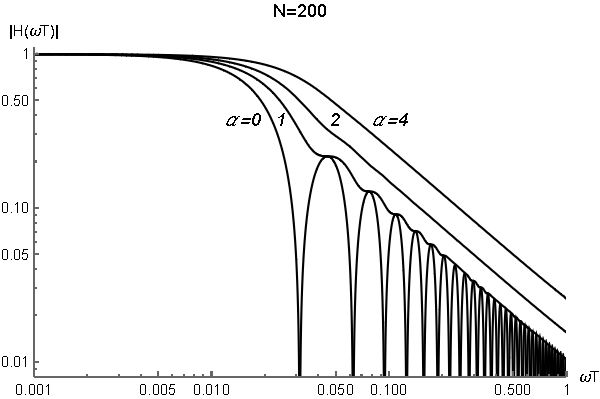} 
\caption{Absolute value of the transfer functions for the orthogonal Hahn filter with $N=200$ and $\alpha=0$,\ $\alpha=1$,\ $\alpha=2$ and $\alpha=4$.}
\end{figure}

\

Taking $\alpha=0$,\ $\alpha=1$ and $\alpha=2$ we get the following $z$-transform functions from \eqref{4.1}
\begin{align*}
&\alpha=0 \qquad H(z,N,0)=\dfrac{z^{1-N}\big(z^N-1\big)}{(N+1)(z-1)} \\
&\alpha=1 \qquad H(z,N,1)=\dfrac{2z^{1-N}\big(z^{N+1}-(N+1)z+N)\big)}{(N+1)(N+2)(z-1)^2} \\
&\alpha=2 \qquad H(z,N,2)=\dfrac{3z^{1-N}\big(2z^{N+2}-(N+1)(N+2)z^2+2N(N+2)z-N(N+1)\big)}{(N+1)(N+2)(N+3)(z-1)^3}
\end{align*}

\

The effect of the number $N$ for the transfer function is given in Figure 4. It is seen that the bandwidth is inversely proportional to the number $N$. The first peak for all functions lies by $\omega\, T=2\pi$. The absolute values of the moduli of the transfer function for $\omega\, T=2\pi$ is the same as the absolute value for $\omega=0$ and is given by \eqref{4.2a}. So for the operation of this filter it must hold that $\omega\, T<1$ provided that the number of points $N$ is greater than about $50$. 
\begin{figure}[ht]
\centering     
\includegraphics[width=10cm]{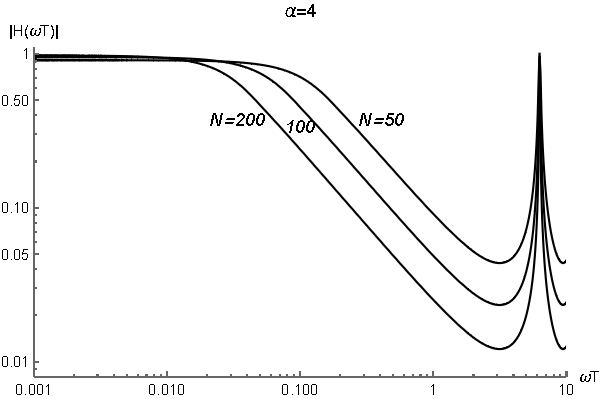} 
\caption{Absolute value of the transfer functions for the orthogonal Hahn filter with $\alpha=4$ and $N=50$,\ $N=100$ and $N=200$.}
\end{figure}

\

Now we have to show that the filter is unbiased. When using the criterion given by \cite[(27)]{10}
\[
\int_0^{2\pi}H(e^{j\omega T})d(\omega T)=\int_0^{2\pi}|H(e^{j\omega T})|^2d(\omega T)
\]
both integrals are zero. This can be proved by the substitution of $z=e^{j\omega\, T}$. Integrating to $z$ gives a result containing only powers of $z$. Back transformation gives the value of $e^{j\omega\, T}$ for the boundaries $\omega=0$ and $\omega=2\pi$ being equal. Then the integrals in the endpoints has the same values. So both integrals are equal to zero. 

\

\

$\text{\bf\Large Acknowledgments}$

\

I thank Prof. Y.S. Shmaliy for his comments on the concept of this paper. 

\

\

\begin{appendices}
$\text{\bf\Large Appendices}$
\section{Transformations of the bounded $_3F_2$ hypergeometric functions with unit argument}
There are seven standard transformations of the bounded $_3F_2$ hypergeometric functions with unit argument \cite[7.4.4. (81-87)]{8}. They are
\begin{align}
\hyp32{-n,a,b}{c,d}{1}
&=\dfrac{(c-a)_n(d-a)_n}{(c)_n(d)_n}\hyp32{-n,a,a+b-c-d-n+1}{a-c-n+1,a-d-n+1}{1} \label{A1} \\
&=\dfrac{(a)_n(c+d-a-b)_n}{(c)_n}(d)_n\hyp32{-n,c-a,d-a}{1-a-n,c+d-a-b}{1} \label{A2} \\
&=\dfrac{(c+d-a-b)_n}{(c)_n}\hyp32{-n,d-a,d-b}{d,c+d-a-b}{1} \label{A3} \\
&=(-1)^n\dfrac{(a)_n(b)_n}{(c)_n(d)_n}\hyp32{-n,1-c-n,1-d-n}{1-a-n,1-b-n}{1} \label{A4} \\
&=(-1)^n\dfrac{(d-a)_n(d-b)_n}{(c)_n(d)_n}\hyp32{-n,1-d-n,a+b-c-d-n}{a-d-n+1,b-d-n+1}{1} \label{A5} \\
&=\dfrac{(c-a)_n}{(c)_n}\hyp32{-n,a,d-b}{d,a-c-n+1}{1} \label{A6} \\
&=\dfrac{(c-a)_n(b)_n}{(c)_n(d)_n}\hyp32{-n,d-b,1-c-n}{1-b-n,a-c-n+1}{1} \label{A7}
\end{align}

\section{Transformations of the Hahn polynomials}
Application of the transformations of the $_3F_2$ hypergeometric functions in Appendix A gives a total of 32 forms of the Hahn polynomials with a lot of conditions. These are
\begin{align*}
Q_{n}(x;\alpha ,\beta ,N)
	&=Q_{n}(-n-\alpha -\beta-1;\alpha ,-x-\alpha -1-n,N) \\
	&=Q_{n}(x;-N-1,n+\alpha +\beta +1+N-n,-\alpha-1) \\
	&=Q_{n}(-n-\alpha -\beta -1;-N-1,-x+N-n,-\alpha-1) \\
	&=\dfrac{(-n-\beta)_n(-N-n-\alpha-\beta -1)_n}{(\alpha +1)_n(-N)_n} \\
	&\qquad\qquad Q_{n}(x-N-\beta ;\beta ,\alpha ,-N-\alpha-\beta -2)
	\qquad\qquad\qquad\qquad\qquad\qquad\qquad\qquad   
\end{align*}
\begin{align*}
Q_{n}(x;\alpha ,\beta ,N)	
	&=\dfrac{(-n-\beta)_n(-N-n-\alpha-\beta -1)_n}{(\alpha +1)_n(-N)_n} \\
	&\qquad\qquad Q_{n}(-n-\alpha -\beta -1;\beta ,N-x-1-n,-N-\alpha -\beta -2) \\
	&=\dfrac{(-n-\beta)_n(-N-n-\alpha-\beta -1)_n}{(\alpha +1)_n(-N)_n} \\
	&\qquad\qquad Q_{n}(-N+x-\beta ;N+\alpha +\beta +1,-N-1,-\beta -1) \\
	&=\dfrac{(-n-\beta)_n(-N-n-\alpha-\beta -1)_n}{(\alpha +1)_n(-N)_n}  \\
	&\qquad\qquad Q_{n}(-N+x-\beta ;N+\alpha +\beta +1,-N-1,-\beta -1)  \\
	&=\dfrac{(-n-\beta)_n(-N-n-\alpha -\beta -1)_n}{(\alpha +1)_n(-N)_n} \\
	&\qquad\qquad Q_{n}(-n-\alpha -\beta -1;N+\alpha +\beta+1,-n-x-\alpha -2,-\beta -1) \\
	&=\dfrac{(n+\alpha +\beta +1)_n(x-n-N-\beta)_n}{(\alpha +1)_n(-N)_n} \\
	&\qquad\qquad Q_{n}(N+n+\alpha +\beta +1;-2n-\alpha -\beta -1,\alpha,N+n-x+\beta) \\
	&=\dfrac{(n+\alpha +\beta +1)_n(x-n-N-\beta)_n}{(\alpha +1)_n(	-N)_n} \\
	&\qquad\qquad Q_{n}(n+\beta ;-2n-\alpha -\beta -1,-N-1,N+n-x+\beta) \\
	&=\dfrac{(n+\alpha +\beta +1)_n(x-n-N-\beta)_n}{(\alpha +1)_n(-N)_n} \\
	&\qquad\qquad Q_{n}(N+n+\alpha +\beta +1;x-n-N-\beta-1,N-n-x,2n+\alpha +\beta)\\
	&=\dfrac{(n+\alpha +\beta +1)_n(x-n-N-\beta)_n}{(\alpha +1)_n(	-N)_n} \\
	&\qquad\qquad Q_{n}(n+\beta ;x-n-N-\beta -1,-n-x-\alpha -1,2n+\alpha+\beta) \\
	&=\dfrac{(x-n-N-\beta)_n}{(\alpha +1)_n}
	Q_{n}(N-x;-N-1,-2n-\alpha-\beta-1,N+n-x+\beta) \\
	&=\dfrac{(x-n-N-\beta)_n}{(\alpha +1)_n}
	Q_{n}(N+n+\alpha +\beta+1;-N-1,x-n,N+n-x+\beta) \\
	&=\dfrac{(x-n-N-\beta)_n}{(\alpha +1)_n}Q_{n}(N-x;-N-n+x-\beta -1,-n-x-\alpha-1,N) \\
	&=\dfrac{(x-n-N-\beta)_n}{(\alpha +1)_n}
	Q_{n}(N+n+\alpha +\beta +1;-N-n+x-\beta-1,\beta ,N) \\
	&=(-1) ^{n}\dfrac{(n+\alpha +\beta+1)_n(-x)_n}{(\alpha +1)_n(-N)_n}
	Q_{n}(n-N-1;-2n-\alpha -\beta -1,\beta ,n-x-1) \\
	&=(-1) ^{n}\dfrac{(n+\alpha +\beta+1)_n(-x)_n}{(\alpha +1)_n(-N)_n}
	Q_{n}(n+\alpha ;-2n-\alpha -\beta -1,N+\alpha +\beta +1,n-x-1) \\
	&=(-1) ^{n}\dfrac{(n+\alpha +\beta+1)_n(-x)_n}{(\alpha +1)_n(-N)_n}Q_{n}(-N+n-1;x-n,-n-x-\alpha -1,2n+\alpha +\beta) \\
	&=(-1) ^{n}\dfrac{(n+\alpha +\beta+1)_n(-x)_n}{(\alpha +1)_n(-N)_n}
	Q_{n}(n+\alpha ;x-n,N-n-x,2n+\alpha +\beta)
\end{align*}
\begin{align*}
Q_{n}(x;\alpha ,\beta ,N)	
	&=(-1) ^{n}\dfrac{(-N-n-\alpha -\beta	-1)_n(x-N)_n}{(\alpha +1)_n(-N)_n} \\
	&\qquad\qquad Q_{n}(x-N-\beta;N+\alpha +\beta+1,-2n-\alpha -\beta -1,n-N+x-1) \\
	&=(-1) ^{n}\dfrac{(-N-n-\alpha -\beta	-1)_n(x-N)_n}{(\alpha +1)_n(-N)_n} \\
	&\qquad\qquad Q_{n}(-N+n-1;N+\alpha +\beta +1,-n-x-\alpha-2,n-N+x-1) \\
	&=(-1) ^{n}\dfrac{(-N-n-\alpha -\beta	-1)_n(x-N)_n}{(\alpha +1)_n(-N)_n} \\
	&\qquad\qquad Q_{n}(-N+x-\beta ;N-n-x,x-n,-N-\alpha-\beta -2) \\
	&=(-1) ^{n}\dfrac{(-N-n-\alpha -\beta	-1)_n(x-N)_n}{(\alpha +1)_n(-N)_n} \\
	&\qquad\qquad Q_{n}(-N+n-1;N-n-x,\beta -1,-N-\alpha -\beta-2) \\
	&=\dfrac{(-n-\beta)_n}{(\alpha+1)_n}Q_{n}(N-x;-N-1,\alpha +\beta +1+N,-\beta -1) \\
	&=\dfrac{(-n-\beta)_n}{(\alpha+1)_n}Q_{n}(-n-\alpha -\beta -1;-N-1,x-n,-\beta -1) \\
	&=\dfrac{(-n-\beta)_n}{(\alpha+1)_n}Q_{n}(N-x;\beta ,\alpha ,N) \\
	&=\dfrac{(-n-\beta)_n}{(\alpha+1)_n}
	Q_{n}(-n-\alpha -\beta -1;\beta ,x-n-N-\beta -1,N) \\
	&=\dfrac{(-n-\beta)_n(-x)_{n}}{(\alpha +1)_n(-N)_n}
	Q_{n}(n+\alpha;x-n,-N-1,-\beta -1) \\
	&=\dfrac{(-n-\beta)_n(-x)_{n}}{(\alpha +1)_n(-N)_n}
	Q_{n}(N-x;x-n,-n-x-\alpha -1,-\beta -1) \\
	&=\dfrac{(-n-\beta)_n(-x)_{n}}{(\alpha +1)_n(-N)_n}
	Q_{n}(n+\alpha;\beta , x-n-N-\beta -1,n-x-1) \\
	&=\dfrac{(-n-\beta)_n(-x)_{n}}{(\alpha +1)_n(-N)_n}
	Q_{n}(N-x;\beta ,-2n-\alpha -\beta -1,n-x-1) 
	\qquad\qquad\qquad\qquad  
\end{align*}

\end{appendices}

\

\end{document}